\input amstex
\documentstyle{amsppt}
\pagewidth{5.4in}
\pageheight{7.6in}
\magnification=1200
\TagsOnRight
\NoRunningHeads
\topmatter
\title
\bf Global and touchdown behaviour of the generalized MEMS device equation 
\endtitle
\author
Kin Ming Hui
\endauthor
\affil
Institute of Mathematics, Academia Sinica,\\
Nankang, Taipei, 11529, Taiwan, R. O. C.
\endaffil
\date
Aug 1, 2008
\enddate
\address
e-mail address: kmhui\@gate.sinica.edu.tw
\endaddress
\abstract
We will prove the local and global existence of solutions of the 
generalized micro-electromechanical system (MEMS) equation $u_t
=\Delta u+\lambda f(x)/g(u)$, $u<1$, in $\Omega\times (0,\infty)$, 
$u(x,t)=0$ on $\partial\Omega\times (0,\infty)$, $u(x,0)=u_0$ in $\Omega$,  
where $\Omega\subset\Bbb{R}^n$ is a bounded domain, $\lambda >0$ is a 
constant, $0\le f\in C^{\alpha}(\overline{\Omega})$, $f\not\equiv 0$, 
for some constant $0<\alpha<1$, $0<g\in C^2((-\infty,1))$ such that 
$g'(s)\le 0$ for any $s<1$ and $u_0\in L^1(\Omega)$ with $u_0\le a<1$ 
for some constant $a$. We prove that there 
exists a constant $\lambda^{\ast}=\lambda^{\ast}(\Omega, f,g)>0$ such that 
the associated stationary problem has a solution for any $0\le\lambda
<\lambda^*$ and has no solution for any $\lambda>\lambda^*$. We obtain 
comparison theorems for the generalized MEMS equation. Under a mild 
assumption on the initial value we prove the convergence of global solutions 
to the solution of the corresponding stationary elliptic equation as 
$t\to\infty$ for any $0\le\lambda<\lambda^*$. We also obtain various 
conditions for the existence of a touchdown time $T>0$ for the solution $u$. 
That is a time $T>0$ such that  $\lim_{t\nearrow T}\sup_{\Omega}u(\cdot,t)=1$. 
\endabstract
\endtopmatter
\NoBlackBoxes
\define \pd#1#2{\frac{\partial #1}{\partial #2}}
\define \1{\partial}
\define \2{\overline}
\define \3{\varepsilon}
\define \4{\widetilde}
\define \5{\underline}
\document

Micro-electromechanical systems (MEMS) are widely used nowadays in many 
electronic devices including accelerometers for airbag deployment in cars, 
inkjet printer heads, and the device for the protection of hard disk, etc. 
Interested readers can read the book, Modeling MEMS and NEMS \cite{PB}, 
by J.A.Pelesko and D.H.~Berstein for the mathematical modeling and various 
applications of 
MEMS devices. Due to the importance of MEMS devices it is important to get a
detail analysis of the mathematical models of MEMS devices. In recent years 
there is a lot of study on the evolution and stationary equations arising 
from MEMS devices by P.~Esposito, N.~Ghoussoub, Y.~Guo, Z.~Pan and M.J.~Ward 
\cite{EGhG},\cite{GhG1},\cite{GhG2},\cite{GPW},\cite{G}, N.I.~Kavallaris, 
T.~Miyasita and T.~Suzuki \cite{KMS}, F.~Lin and Y.~Yang \cite{LY},
L.~Ma and J.C.~Wei \cite{MW} and J.A.Pelesko \cite{P}, etc.

Let $\Omega\subset\Bbb{R}^n$ be a bounded $C^2$ domain. Let 
$$
0\le f\in C^{\alpha}(\2{\Omega})\quad\text{ for some constant }0<\alpha<1
\quad\text{ and }f\not\equiv 0\quad\text{ in }\Omega\tag 0.1
$$ 
and let
$$
0<g\in C^2((-\infty,1))\quad\text{ such that }
g'(s)\le 0\quad\forall s<1.\tag 0.2
$$ 
In this paper we will study the generalized MEMS equation
$$\left\{\aligned
&u_t=\Delta u+\frac{\lambda f(x)}{g(u)}\quad\text{ in }\Omega\times (0,T)\\
&u(x,t)=0\qquad\qquad\text{ on }\1\Omega\times (0,T)\\
&u(x,0)=u_0\qquad\quad\text{ in }\Omega
\endaligned\right.\tag 0.3
$$
and the associated stationary problem,
$$\left\{\aligned
&-\Delta v=\frac{\lambda f(x)}{g(v)}\quad\text{ in }\Omega\\
&v(x)=0\qquad\qquad\text{ on }\1\Omega.
\endaligned\right.\tag $S_{\lambda}$
$$
When $g(u)=(1-u)^2$, (0.3) and ($S_{\lambda}$) reduces to the evolution 
and stationary
MEMS equations respectively which were studied extensively in 
\cite{EGhG},\cite{GhG1},\cite{GhG2},\cite{GPW},\cite{G},\cite{P}. 
An equation similar to ($S_{\lambda}$) arising from the motion of thin
films of viscous fluid is studied by H.~Jiang and W.M.~Ni in \cite{JN}. 
The aymptotic and touchdown behaviour of solutions of ($S_{\lambda}$) with
$g(u)=(1-u)^2$ and $u_0\equiv 0$ was studied in \cite{GhG2} and \cite{G}.
When $g(u)=(1-u)^p$ with $p>0$, ($S_{\lambda}$) was studied 
by L.~Ma, J.C.~Wei, Z.~Wang and L.~Ruan \cite{MW},\cite{WR}. 
The equation (0.3) and ($S_{\lambda}$) with 
$g(u)=(1-u)^p$ and $u_0\in [0,1)$ were also studied by N.I.~Kavallaris, 
T.~Miyasita, T.~Suzuki \cite{KMS}. By the results of \cite{GhG1},\cite{GhG2}, 
and \cite{WR}, when $g(u)=(1-u)^p$ with $p>0$, there exists a constant 
$\lambda^*>0$ such that ($S_{\lambda}$) has a solution for any 
$0\le\lambda<\lambda^{\ast}$ and ($S_{\lambda}$) 
has no solution for any $\lambda>\lambda^{\ast}$. 

In this paper we will show that there exists a constant $\lambda^{\ast}>0$ 
such that similar results hold for ($S_{\lambda}$).
The constant $\lambda^{\ast}$ is called the pull-in voltage of the equation
($S_{\lambda}$) in the literature of MEMS.
For any $u_0\in L^1(\Omega)$ with $u_0\le a<1$ for 
some constant $a$ we will prove the local existence and comparison theorems
of solutions of (0.3). If $u$ is a global solution of (0.3) with 
$0\le\lambda<\lambda^{\ast}$, then under a mild
assumption on the initial value we prove the convergence of the solution 
of (0.3) as $t\to\infty$. We also obtain various conditions for the solution 
$u$ of (0.3) to touchdown at a finite time. That is the existence of a 
time $T>0$ such that
$$
\lim_{t\nearrow T}\sup_{\Omega}u(\cdot,t)=1.
$$
The plan of the paper is as follows. In section 1 we will prove the 
existence of finite pull-in voltage $\lambda^{\ast}>0$ of ($S_{\lambda}$)
and the existence and non-existence of solutions of ($S_{\lambda}$). 
We will also prove the non-existence of bounded solution of the 
stationary problem in $\Bbb{R}^n$. In section 2
we will prove the existence of solutions and various comparsion results
for solutions of (0.3). In section 3 we will prove the global convergence 
of solutions of (0.3) for $0\le\lambda<\lambda^{\ast}$. We also obtain 
various conditions for the solutions of (0.3) to have finite touchdown time. 

We start with a definition. We say that $v$ is a solution (subsolution,
supersolution respectively) of ($S_{\lambda}$) if $v\in C^2(\Omega)\cap 
C(\2{\Omega})$, $v<1$ in $\Omega$, satisfies 
$$
-\Delta v=\frac{\lambda f(x)}{g(v)}\quad\text{ in }\Omega
$$
($\le$, $\ge$ respectively) with $v(x)=0$ ($\le$, $\ge$ respectively) on 
$\1\Omega$. Note that by the maximum principle for superharmonic function
if $v$ is a solution or supersolution of ($S_{\lambda}$), then
$v\ge 0$ in $\2{\Omega}$. We say that $v$ is a minimal
solution of ($S_{\lambda}$) if $v$ is a solution of ($S_{\lambda}$)
and $v\le\4{v}$ in $\Omega$ for any solution $\4{v}$ of ($S_{\lambda}$). 

For any 
$$
u_0\in L^1(\Omega)\text{ with }u_0\le a\text{ on }
\Omega\tag 0.4
$$
for some $a\in (0,1)$
we say that $u$ is a solution (subsolution, supersolution 
respectively) of (0.3) in $\Omega\times (0,T)$ if $u\in C^{2,1}
(\Omega\times (0,T))\cap C(\2{\Omega}\times (0,T))$ satisfies 
$$
u_t=\Delta u+\frac{\lambda f(x)}{g(u)}\quad\text{ in }\Omega\times (0,T)
$$
($\le$, $\ge$ respectively) with $u(x,t)=0$ ($\le$, $\ge$ respectively) on 
$\1\Omega\times (0,T)$, 
$$
\sup_{\2{\Omega}\times (0,T']}u(x,t)<1\quad\forall 0<T'<T
$$ 
and 
$$
\|u(\cdot,t)-u_0\|_{L^1(\Omega)}\to 0\quad\text{ as }t\to 0.\tag 0.5 
$$
For any solution $u$ of (0.3) we define the touchdown time
$T_{\lambda}=T_{\lambda}(\Omega,f,g)>0$ as the time which satisfies 
$$\left\{\aligned
\sup_{\Omega}u(x,t)&<1\quad\forall 0<t<T_{\lambda}\\
\lim_{t\nearrow T_{\lambda}}\sup_{\Omega}u(x,t)&=1.
\endaligned\right.
$$
We say that $u$ has a finite touchdown time if $T_{\lambda}<\infty$ and we
say that $u$ touchdowns at time infinity if $T_{\lambda}=\infty$.

Let $G(x,y,t)$, $x,y\in\Omega$, $t>0$, be the Dirichlet Green function of 
the heat equation in $\Omega\times (0,\infty)$. That is for any $y\in\Omega$,
$$\left\{\aligned 
&\1_tG=\Delta_xG\quad\text{ in }\Omega\times (0,\infty)\\
&G(x,y,t)=0\quad\forall x\in\1\Omega, t>0\\
&\lim_{t\to 0}G(x,y,t)=\delta_y
\endaligned\right.
$$
where $\delta_y$ is the delta mass at $y$. By the maximum principle,
$$
0\le G(x,y,t)\le\frac{1}{(4\pi)^{\frac{n}{2}}}e^{-|x-y|^2/4t}.\tag 0.6
$$
For any $K\subset R^n\times (0,\infty)$, $0<\beta<1$, let
$$
C^{2,1}(K)=\{f: f,f_t, f_{x_i}, f_{x_i,x_j}\in C(K)\quad\forall 
i,j=1,2,\dots,n\}
$$
and let $C^{2+\beta,1+(\beta/2)}(K)$ denote the class of
all functions $f\in C^{2,1}(K)$ such that
$$\left\{\aligned
&|f_{x_i,x_j}(x_1',t_1')- f_{x_i,x_j}(x_2',t_2')|
\le C(|x_1'-x_2'|^{\beta}+|t_2'-t_1'|^{\beta/2})\}\quad\forall
(x_1',t_1'), (x_2',t_2')\in K\\
&|f_t(x_1',t_1')- f_t(x_2',t_2')|
\le C(|x_1'-x_2'|^{\beta}+|t_2'-t_1'|^{\beta/2})\}\qquad\qquad\forall
(x_1',t_1'), (x_2',t_2')\in K\endaligned\right. 
$$
holds for some constant $C>0$ and any $i,j=1,2,\cdots ,n$.

For any set $A$, let $\chi_A$ be the characteristic function of $A$. 
For any $a\in\Bbb{R}$, let $a_-=\max (0,-a)$. For any $x_0\in\Bbb{R}^n$,
$R>0$, let $B_R(x_0)=\{x\in\Bbb{R}^n:|x-x_0|<R\}$ and $B_R=B_R(0)$.
Let $\Cal{C}$ be the family of bounded $C^2$ domain $\Omega_1\subset
\Bbb{R}^n$ such that $\2{\Omega}\subset\Omega_1$. For any $\Omega_1\in
\Cal{C}$ let $\mu_{\Omega_1}>0$ be the first Dirichlet 
eigenvalue of $-\Delta$ in $\Omega_1$ and $\psi_{\Omega_1}$ be the 
corresponding positive eigenfunction normalized such that 
$$
\max_{\Omega_1}\psi_{\Omega_1}=1\quad\text{ and }\quad s_{\Omega_1}
=\min_{\2{\Omega}}\psi_{\Omega_1}>0.
$$ 
Let $\mu_1>0$ be the first Dirichlet eigenvalue of $-\Delta$ in 
$\Omega$ and let $\phi_1$ be the first positive Dirichlet 
eigenfunction of $-\Delta$ in $\Omega$ normalized such that 
$\int_{\Omega}\phi_1\,dx=1$. Let
$$
\nu_{\Omega}=\sup_{\Omega_1\in\Cal{C}}\mu_{\Omega_1}s_{\Omega_1}.
$$

$$
\text{\bf Section 1}
$$

In this section we will prove the existence of finite pull-in voltage 
$\lambda^{\ast}>0$ of ($S_{\lambda}$) and the existence and non-existence of 
solutions of ($S_{\lambda}$). We also obtain various estimates for
$\lambda^{\ast}$. 

\proclaim{Theorem 1.1}
Suppose $f$ satisfies (0.1) and $g$ satisfies (0.2). Then
there exists a constant $\lambda^{\ast}=\lambda^{\ast}(\Omega, f,g)>0$ such 
that 
$$\aligned
&\text{ (i) }\forall 0\le\lambda<\lambda^{\ast},\text{ there exists at least  
one solution of ($S_{\lambda}$)}\\
&\text{ (ii) }\forall\lambda>\lambda^{\ast},\text{ there exists no solution of 
($S_{\lambda}$)}.\endaligned
$$
Moreover
$$
\nu_{\Omega}\frac{\sup_{0<s<1}sg(s)}{\max_{\2{\Omega}}f}
\le\lambda^{\ast}\le\frac{\mu_1g(0)}{\int_{\Omega}f\phi_1\,dx}.
\tag 1.1
$$
\endproclaim
\demo{Proof}
Since the proof of the theorem is similar to the proof of Theorem 2.1 of 
\cite{GhG1}, we will sketch the argument here. Note that $v\equiv 0$ in
$\Omega$ is a solution of ($S_{\lambda}$) when $\lambda=0$. Let 
$D=\{\lambda>0:(S_{\lambda})\text{ has a solution}\}$ and 
$$
\lambda^{\ast}=\lambda^{\ast}(\Omega, f,g)=\sup_{\lambda\in D}\lambda.
$$ 
We claim that $D\ne\phi$. In order to prove the claim we first observe that 
$v\equiv 0$ on $\2{\Omega}$ is a subsolution of $(S_{\lambda})$ for any 
$\lambda>0$. We will next construct a supersolution
of $(S_{\lambda})$. For any $\Omega_1\in\Cal{C}$ and $0<A<1$ let 
$\psi=A\psi_{\Omega_1}$. Then by (0.2) for any 
$$
0\le\lambda\le\mu_{\Omega_1}s_{\Omega_1}\frac{Ag(A)}{\max_{\2{\Omega}}f},
$$
we have
$$
-\Delta\psi=A\mu_{\Omega_1}\psi_{\Omega_1}\ge
 \frac{\lambda f}{g(A\psi_{\Omega_1})}\quad\text{ in }\Omega.
$$
Hence $\psi$ is a supersolution of ($S_{\lambda}$). Let $v_0\equiv 0$
in $\Omega$ and for any $k\ge 1$, let $v_k$ be the solution of
$$\left\{\aligned
-\Delta v_k&=\frac{\lambda f(x)}{g(v_{k-1})}\quad\text{ in }\Omega\\
v_k(x)&=0\qquad\qquad\text{on }\1\Omega.
\endaligned\right.
$$
By (0.2) and an argument similar to that of \cite{GhG1}, 
$0\le v_k\le v_{k+1}\le\psi <1$ in $\Omega$ for all $k\ge 0$ 
and $v_k$ will converge to the minimal solution $v$ of ($S_{\lambda}$)
as $k\to\infty$. Hence $D\ne\phi$ and the left hand side inequality of 
(1.1) holds. 

Suppose now $v$ is a solution ($S_{\lambda}$). Multiplying ($S_{\lambda}$) 
by $\phi_1$ and integrating over $\Omega$, by (0.2) we have
$$
\mu_1\ge\mu_1\int_{\Omega}v\phi_1\,dx
=-\int_{\Omega}v\Delta\phi_1\,dx
=-\int_{\Omega}\phi_1\Delta v\,dx
=\lambda\int_{\Omega}\frac{f\phi_1}{g(v)}\,dx
\ge\frac{\lambda}{g(0)}\int_{\Omega}f\phi_1\,dx.
$$
Hence
$$
\lambda\le\frac{\mu_1g(0)}{\int_{\Omega}f\phi_1\,dx}.
$$
Thus the right hand side inequality of (1.1) and (ii) follows. For any 
$0\le\lambda<\lambda^{\ast}$, there exists 
$\lambda<\lambda_1<\lambda^{\ast}$ such that ($S_{\lambda_1}$) has a solution 
$v_{\lambda_1}$. Then $v_{\lambda_1}$ is a supersolution of ($S_{\lambda}$).
By (0.2) and the monotone iteration scheme as before (cf. \cite{GhG1}) 
($S_{\lambda}$) has a solution $v$ satisfying $0\le v\le v_{\lambda}$ in 
$\Omega$ and (i) follows.
\enddemo

We will now let $\lambda^{\ast}$ be given by Theorem 1.1 for the rest of the 
paper. The following result improves the upper bound of $\lambda^{\ast}$
of Theorem 1.1.

\proclaim{Proposition 1.2}
Suppose $f$ satisfies (0.1) and $g$ satisfies (0.2). Then
$$
\lambda^{\ast}\le\mu_1
\frac{\int_0^1g(s)\,ds}{\int_{\Omega}f\phi_1\,dx}\tag 1.2
$$
where
$$
H(v)=\int_v^1g(s)\,ds.\tag 1.3
$$
\endproclaim
\demo{Proof}
Suppose $v$ is a solution of ($S_{\lambda}$). 
Multiplying ($S_{\lambda}$) by $g(v)\phi_1$ and integrating over $\Omega$,
$$\align
\lambda\int_{\Omega}f\phi_1\,dx=&-\int_{\Omega}g(v)\phi_1\Delta v\,dx\\
=&\int_{\Omega}\nabla (g(v)\phi_1)\cdot\nabla v\,dx
-\int_{\1\Omega}g(v)\phi_1\frac{\1 v}{\1\nu}\,d\sigma\\
=&\int_{\Omega}g'(v)\phi_1|\nabla v|^2\,dx
+\int_{\Omega}g(v)\nabla\phi_1\cdot\nabla v\,dx\\
\le&-\int_{\Omega}\nabla\phi_1\cdot\nabla H(v)\,dx\\
=&\int_{\Omega}H(v)\Delta\phi_1\,dx-\int_{\1\Omega}H(v)\frac{\1\phi_1}{\1\nu}
\,d\sigma\\
=&-\mu_1\int_{\Omega}H(v)\phi_1\,dx-H(0)\int_{\1\Omega}
\frac{\1\phi_1}{\1\nu}\,d\sigma\\
=&-\mu_1\int_{\Omega}H(v)\phi_1\,dx-H(0)\int_{\Omega}\Delta\phi_1
\,dx\\
=&-\mu_1\int_{\Omega}H(v)\phi_1\,dx+\mu_1H(0)
\int_{\Omega}\phi_1\,dx\\
=&-\mu_1\int_{\Omega}H(v)\phi_1\,dx+\mu_1H(0)\\
\endalign
$$
and (1.2) follows.  
\enddemo

We will next prove a more computable bound for $\lambda^{\ast}$.

\proclaim{Proposition 1.3}
Suppose $f\in C^1(\2{\Omega})$ satisfies 
$$
\delta_1=\inf_{\Omega}f>0,\tag 1.4
$$
$g$ satisfies (0.2), and $\Omega\subset\Bbb{R}^n$, $n\ge 2$, 
is a strictly star-shape domain 
such that $x\cdot\nu\ge b>0$ on $\1\Omega$ where $\nu$ is the unit outward 
normal to $\1\Omega$ at $x\in\1\Omega$. Then
$$
\lambda^{\ast}\le\frac{((n+2)\|f\|_{L^{\infty}}+2b_1)
|\1\Omega|}{\delta_1^2b|\Omega|}
g(0)
$$
where $b_1=\sup_{\Omega}|x\cdot\nabla f|$. In particular if $\Omega=B_R$, 
then
$$
\lambda^{\ast}\le\frac{n((n+2)\|f\|_{L^{\infty}}+2b_1)}{\delta_1^2R}g(0).
$$ 
\endproclaim
\demo{Proof}
Suppose $\lambda>0$ and $v$ is a solution of ($S_{\lambda}$). 
By ($S_{\lambda}$) and the Pohozaev identity \cite{N},
$$\align
&\frac{1}{2}\int_{\1\Omega}(x\cdot\nu)\left(\frac{\1 v}{\1\nu}\right)^2
\,d\sigma\\
=&\lambda n\int_{\Omega}f(x)\left(\int_0^{v(x)}\frac{ds}{g(s)}\right)\,dx
-\lambda\frac{(n-2)}{2}\int_{\Omega}\frac{vf(x)}{g(v)}\,dx\\
&\qquad +\lambda\int_{\Omega}(x\cdot\nabla f(x))\left(\int_0^{v(x)}
\frac{ds}{g(s)}\right)\,dx.\tag 1.5
\endalign
$$
By (0.2),
$$
\left(\int_0^v\frac{ds}{g(s)}\right)\le\frac{v}{g(v)}.
$$
Hence the right hand side of (1.5) is less than
$$\align
\le&\lambda n\int_{\Omega}\frac{vf(x)}{g(v)}\,dx
-\lambda\frac{(n-2)}{2}\int_{\Omega}\frac{vf(x)}{g(v)}\,dx
+\lambda b_1\int_{\Omega}\frac{v}{g(v)}\,dx\\
\le&\lambda\left(\frac{(n+2)}{2}\|f\|_{L^{\infty}}+b_1\right)
\int_{\Omega}\frac{dx}{g(v)}.\tag 1.6
\endalign
$$
Now by the Holder inequality, the Green theorem and ($S_{\lambda}$),
$$
\int_{\1\Omega}(x\cdot\nu)\left(\frac{\1 v}{\1\nu}\right)^2\,d\sigma
\ge\frac{b}{|\1\Omega|}\left(\int_{\1\Omega}\frac{\1 v}{\1\nu}\,d\sigma
\right)^2
=\frac{b}{|\1\Omega|}\left(\int_{\Omega}\Delta v\,dx\right)^2
\ge\frac{b\lambda^2\delta_1^2}{|\1\Omega|}\left(\int_{\Omega}\frac{dx}{g(v)}
\right)^2.\tag 1.7
$$
By (1.5), (1.6) and (1.7),
$$\align
\frac{b\lambda^2\delta_1^2}{2|\1\Omega|}\left(\int_{\Omega}\frac{dx}{g(v)}
\right)^2
\le&\lambda\left(\frac{(n+2)}{2}\|f\|_{L^{\infty}}+b_1\right)
\int_{\Omega}\frac{dx}{g(v)}\\
\Rightarrow\quad\frac{((n+2)\|f\|_{L^{\infty}}+2b_1)}{b}|\1\Omega|
\ge&\lambda\delta_1^2\int_{\Omega}\frac{dx}{g(v)}
\ge\lambda\delta_1^2\frac{|\Omega|}{g(0)}\\
\Rightarrow\qquad\qquad\qquad\qquad\qquad\qquad\lambda
\le&\frac{((n+2)\|f\|_{L^{\infty}}+2b_1)|\1\Omega|}{\delta_1^2b|\Omega|}g(0)\\
\Rightarrow\qquad\qquad\qquad\qquad\qquad\quad\,\,\,\lambda^{\ast}
\le&\frac{((n+2)\|f\|_{L^{\infty}}+2b_1)|\1\Omega|}{\delta_1^2b|\Omega|}g(0)
\endalign
$$
and the proposition follows.
\enddemo

\proclaim{Corollary 1.4}
Let $f\in C^1(\2{\Omega})$ satisfy (1.4) such that supp$\,\nabla f\subset
B_{R_1}$ for some constant $R_1>1$ and let $g$ satisfy (0.2).
For any $\lambda>0$ there does not exist any bounded solution for the problem,
$$
-\Delta w=\frac{\lambda f(x)}{g(w)}, w<1,\quad\text{ in }\Bbb{R}^n
\tag 1.8
$$
\endproclaim
\demo{Proof}
Suppose there exists $\lambda>0$ such that (1.8) has a bounded solution $w$.
Without loss of generality we may assume that $0\le w<1$ in $\Bbb{R}^n$. 
Let 
$$
R_2=\frac{2n((n+2)\|f\|_{L^{\infty}}+2R_1\|\nabla f\|_{L^{\infty}})}
{\delta_1^2\lambda}g(0).
$$
By Proposition 1.3 $\lambda^{\ast}(B_{R_2},f,g)\le\lambda/2$. On the other
hand since $w$ is a supersolution of ($S_{\lambda}$) with $\Omega=B_{R_2}$,
by the construction of solutions of ($S_{\lambda}$) in Theorem 1.1,
there exists a solution $v$ of ($S_{\lambda}$) with $\Omega=B_{R_2}$ 
satisfying $0\le v\le w$. Hence $\lambda^{\ast}(B_{R_2},f,g)\ge\lambda$
and contradiction arises. Thus no such solution $w$ exists.
\enddemo

\proclaim{Proposition 1.5}
Let $\Omega_1\subset\Omega_2$ and let $f_1$, $f_2$ satisfy (0.1) in
$\Omega_1$, $\Omega_2$, respectively for some constant $0<\alpha<1$
such that $f_1\le f_2$ in $\Omega_1$. Let $g_1,g_2$ satisfy (0.2) such that
$g_1(s)\ge g_2(s)>0$ for any $s<1$. Then $\lambda^{\ast}(\Omega_1,f_1,g_1)
\ge\lambda^{\ast}(\Omega_2,f_2,g_2)$.
If $0\le\lambda<\lambda^{\ast}(\Omega_2,f_2,g_2)$ and $v_1$, $v_2$, are the 
minimal solutions of ($S_{\lambda}$) with $\Omega=\Omega_1,\Omega_2$, 
$f=f_1,f_2, g=g_1, g_2$, respectively, then $v_1\le v_2$ in $\Omega_1$. 
If moreover $\Omega_1=\Omega_2=\Omega$ and $f_1\not\equiv f_2$, 
then $v_1<v_2$ in $\Omega$.
\endproclaim
\demo{Proof}
For any $\lambda<\lambda^{\ast}(\Omega_2,f_2,g_2)$, let $v_2$ be the minimal
solution of ($S_{\lambda}$) with $\Omega=\Omega_2, f=f_2, g=g_2$. Then 
$v_2$ is a
supersolution of ($S_{\lambda}$) with $\Omega=\Omega_1, f=f_1, g=g_1$. Since
$0$ is a subsolution of ($S_{\lambda}$) with $\Omega=\Omega_1, f=f_1, g=g_1$,
by the monotone iteration scheme for the construction of solution of
($S_{\lambda}$) as in the proof of Theorem 1.1 the minimal solution 
$v_1$ of ($S_{\lambda}$) with $\Omega=\Omega_1, f=f_1, g=g_1$ satisfies
$0\le v_1\le v_2$ in $\Omega_1$. Hence $\lambda^{\ast}(\Omega_1,f_1,g_1)
\ge\lambda^{\ast}(\Omega_2,f_2,g_2)$. 

We next suppose that $\Omega_1=\Omega_2=\Omega$ and $f_1\not\equiv f_2$.
Let $G(x,y)$ be the Green function for $\Delta$ in $\Omega$.
Then
$$
v_i(x)=\lambda\int_{\Omega}G(x,y)\frac{f_i(y)}{g_i(v_i(y))}\,dx
\quad\forall i=1,2.\tag 1.9
$$
Since $v_1\le v_2$ in $\Omega$, by (0.2) $f_1(x)/g_1(v_1)\le f_2(x)/g_2(v_2)$ 
in $\Omega$. If $f_1\not\equiv f_2$, there exists a set $A\subset\Omega$ 
of positive measure such that $f_1(x)/g_1(v_1)<f_2(x)/g_2(v_2)$ in $A$. 
Then by (1.9), $v_1<v_2$ in $\Omega$ and the proposition follows.
\enddemo

For any solution $v$ of ($S_{\lambda}$) we let
$$
L_{v,\lambda}w=-\Delta w+\lambda\frac{f(x)g'(v)}{g(v)^2}w
$$
be the linearized operator of ($S_{\lambda}$) around the solution $v$.
Let 
$$
\4{\mu}_1=\4{\mu}_1(\lambda,v)
=\inf_{w\in H_0^1(\Omega)}\frac{\int_{\Omega}|\nabla w|^2\,dx
+\lambda\int_{\Omega}(fg'(v)/g(v)^2)w^2\,dx}{\int_{\Omega}w^2\,dx}
$$
and $\4{\phi}_1$ be the first eigenvalue and the corresponding first 
positive eigenfunction of $L_{v,\lambda}$. We say that $v$ is a stable
solution of ($S_{\lambda}$) if $v$ is a solution of ($S_{\lambda}$)
with $\4{\mu}_1(\lambda,v)>0$.

\proclaim{Theorem 1.6}
Let $f$ satisfy (0.1) and $g$ satisfy (0.2) and
$$
\left(\frac{1}{g}\right)''(s)\ge 0\quad\forall s<1.\tag 1.10
$$
Suppose $v$ and $\4{v}$ are solution and supersolution of
($S_{\lambda}$) respectively. If $\4{\mu}_1=\4{\mu}_1(\lambda,v)>0$, 
then $\4{v}\ge v$ in $\Omega$. If $\4{\mu}_1=0$, then $\4{v}\equiv v$ in 
$\Omega$.  
\endproclaim
\demo{Proof}
We will use a modification of the proof of Lemma 4.1 of \cite{GhG1} to prove
the theorem. Let
$$
h(x,s)=-\Delta (s\4{v}+(1-s)v)-\frac{\lambda f}{g(s\4{v}+(1-s)v)}
\quad\forall 0\le s\le 1.
$$
Then 
$$
h(x,0)=0.\tag 1.11
$$
By (1.10) and the Jensen inequality,
$$
-\Delta (s\4{v}+(1-s)v)=\lambda f\left(\frac{s}{g(\4{v})}
+\frac{1-s}{g(v)}\right)\ge\frac{\lambda f}{g(s\4{v}+(1-s)v)}\quad\text{ in }
\Omega\quad\forall 0\le s\le 1.
$$
Hence
$$
h(x,s)\ge 0\quad\text{ in }\Omega\quad\forall 0\le s\le 1.\tag 1.12
$$
By (1.11) and (1.12),
$$
\frac{\1 h}{\1 s}(x,0)\ge 0\quad
\Rightarrow\quad -\Delta (\4{v}-v)+\lambda f\frac{g'(v)}{g(v)^2}(\4{v}-v)\ge 0
\quad\text{ in }\Omega.\tag 1.13
$$
Suppose first $\4{\mu}_1>0$. Multiplying (1.13) by $(\4{v}-v)_-$ and 
integrating over $\Omega$,
$$\align
0&\ge\int_{\Omega}|\nabla (\4{v}-v)_-|^2\,dx
+\lambda\int_{\Omega}f\frac{g'(v)}{g(v)^2}(\4{v}-v)_-^2\,dx\\
&\ge\4{\mu}_1\int_{\Omega}(\4{v}-v)_-^2\,dx\\
\Rightarrow\quad&\4{v}\ge v\quad\text{ in }\Omega.\tag 1.14
\endalign
$$
Suppose now $\4{\mu}_1=0$. Multiplying (1.13) by $\4{\phi}_1$ and
integrating over $\Omega$,
$$\align
0\le&-\int_{\Omega}\4{\phi}_1\Delta (\4{v}-v)\,dx
+\lambda\int_{\Omega}\4{\phi}_1f\frac{g'(v)}{g(v)^2}(\4{v}-v)\,dx\\
=&\int_{\Omega}(\4{v}-v)\left(-\Delta\4{\phi}_1
+\lambda f\frac{g'(v)}{g(v)^2}\4{\phi}_1\right)\,dx\\
=&0.\tag 1.15
\endalign
$$
Hence by (1.13), (1.15) and the positivity of $\4{\phi}_1$ in $\Omega$,
$$
\frac{\1 h}{\1 s}(x,0)
=-\Delta (\4{v}-v)+\lambda f\frac{g'(v)}{g(v)^2}(\4{v}-v)=0
\quad\text{ in }\Omega.\tag 1.16
$$
By (1.10), (1.11), (1.12) and (1.16),
$$\align
\frac{\1^2h}{\1 s^2}(x,0)\ge 0\quad\Rightarrow\quad&
-\lambda f\left(\frac{1}{g}\right)''(v)(\4{v}-v)^2\ge 0\quad\text{ in }
\Omega\\
\Rightarrow\quad&\4{v}=v\quad\text{ in }\Omega\setminus D_1
\endalign
$$
where $D_1=\{x\in\Omega: f(x)=0\}$. By (1.16) $\Delta (\4{v}-v)=0$
in $D_1$. Since $\4{v}-v=0$ on $\1 D_1$, $\4{v}\equiv v$ on $D_1$. Hence
$\4{v}\equiv v$ in $\Omega$ and the theorem follows.
\enddemo

By Theorem 1.6 and an argument similar to the proof of Theorem 4.2
of \cite{GhG1} we have the following theorem.

\proclaim{Theorem 1.7}
Let $f$ satisfy (0.1) and $g$ satisfy (0.2) and (1.10).
For each $0<\lambda<\lambda^{\ast}$ let $v_{\lambda}$ be the minimal solution
of ($S_{\lambda}$). Then $v_{\lambda}(x)$ is a stable solution 
of ($S_{\lambda}$) for any $0<\lambda<\lambda^{\ast}$. Moreover for 
each $x\in\Omega$, $v_{\lambda}(x)$ 
is differentiable and strictly increasing with respect to $\lambda\in
(0,\lambda^{\ast})$ and $\4{\mu}_1(\lambda,v_{\lambda})$ is a decreasing
function of $\lambda\in (0,\lambda^{\ast})$.
\endproclaim

\proclaim{Proposition 1.8}
Let $f$ satisfy (0.1) and $g$ satisfy (0.2) and (1.10).
For each $0<\lambda<\lambda^{\ast}$ let $v_{\lambda}$ be the minimal solution
of ($S_{\lambda}$). Suppose $v$ is a solution of ($S_{\lambda}$) and 
$v\not\equiv v_{\lambda}$. Then $\4{\mu}_1(\lambda,v)<0$ and the function
$w=v-v_{\lambda}$ is in the negative space of $L_{v,\lambda}$.
\endproclaim
\demo{Proof}
Since $v_{\lambda}$ is the minimal solution of ($S_{\lambda}$),
$v\ge v_{\lambda}$ in $\Omega$. Let $D_1=\{x\in\Omega: f(x)=0\}$
and $D_2=\{x\in \Omega\setminus D_1:v(x)\ne v_{\lambda}(x)\}$. 
If $v\equiv v_{\lambda}$ in $\Omega\setminus D_1$, then 
$\Delta (v-v_{\lambda})=0$ in $D_1$ and $v=v_{\lambda}$ on $\1 D_1$. 
Thus $v\equiv v_{\lambda}$ on $\2{\Omega}$. Contradiction arises. Hence
$v\not\equiv v_{\lambda}$ in $\Omega\setminus D_1$ and $D_2$ is a set of 
positive measure. By the mean value theorem,
$$\align
L_{v,\lambda}(v-v_{\lambda})
=&-\Delta (v-v_{\lambda})-\lambda f\left(\frac{1}{g}\right)'(v)
(v-v_{\lambda})\\
=&\lambda f\left\{\frac{1}{g(v)}-\frac{1}{g(v_{\lambda})}
-\left(\frac{1}{g}\right)'(v)(v-v_{\lambda})\right\}\\
=&\lambda f\left\{\left(\frac{1}{g}\right)'(\xi_1)
-\left(\frac{1}{g}\right)'(v)\right\}(v-v_{\lambda})\\
=&\lambda f\left(\frac{1}{g}\right)''(\xi_2)(v-v_{\lambda})(\xi_1-v)
\quad\text{ in }D_2\tag 1.17
\endalign
$$
for some functions $\xi_1(x)\in (v_{\lambda}(x),v(x))$, $\xi_2(x)\in 
(v_{\lambda}(x),\xi_1(x))$. Hence by (1.10) and (1.17),
$$
<L_{v,\lambda}w,w>
=\int_{D_2}\lambda f\left(\frac{1}{g}\right)''(\xi_2)
(v-v_{\lambda})^2(\xi_1-v)\,dx<0.
$$
Thus $\4{\mu}_1(\lambda,v)<0$ and the proposition follows.
\enddemo

$$
\text{\bf Section 2}
$$

In this section we will prove the local and global existence of solutions 
of (0.3). We also obtain various comparison results for the solutions 
of (0.3).

\proclaim{Theorem 2.1}
Let $u_{0,1}, u_{0,2}\in L^1(\Omega)$. Let $f\in C(\2{\Omega})$ and 
$0<g\in C^2((-\infty,1))$. Suppose $u_1$, $u_2$, are subsolution
and supersolution of (0.3) in $\Omega\times (0,T)$ with initial value $u_0
=u_{0,1}, u_{0,2}$, respectively such that 
$$
a_1=\max(\sup_{\2{\Omega}\times (0,T)}u_1(x,t),
\sup_{\2{\Omega}\times (0,T)}u_2(x,t))<1.\tag 2.1
$$ 
Suppose either (1.10) holds or there exists $a_2<1$ such that
$$
u_1(x,t),u_2(x,t)\ge a_2\quad\text{ on }\Omega\times (0,T).\tag 2.2
$$
Then
$$
\text{(i)}\quad\int_{\Omega}(u_1-u_2)_+(x,t)dx\le e^{bt}
\int_{\Omega}(u_{0,1}-u_{0,2})_+dx\quad\forall 0\le t<T
$$
hold for some constant $b>0$ depending on $\lambda$, $f$, and $a_1$ if 
(1.10) holds and on $\lambda$, $f$, $a_1$ and $a_2$ if (2.2) holds. If 
both $u_1$ and $u_2$ are solutions of (0.3) in $\Omega\times (0,T)$
with initial value $u_0=u_{0,1}, u_{0,2}$, respectively, then
$$  
\text{(ii)}\quad\int_{\Omega}|u_1-u_2|(x,t)dx\le e^{bt}
\int_{\Omega}|u_{0,1}-u_{0,2}|dx\quad\forall 0\le t<T.
$$
\endproclaim
\demo{Proof}
We will use a modification of the technique of Dahlberg and C.~Kenig
\cite{DK} to prove the theorem.
Let $h\in C_0^{\infty}(\Omega)$ be such that $0\le h\le 1$. For any 
$t_1\in (0,T)$, let $\eta$ be the solution of 
$$\left\{\aligned
\eta_t+\Delta\eta+H\eta&=0\qquad\text{ in }\Omega\times (0,t_1)\\
\eta&=0\qquad\text{ on }\1\Omega\times (0,t_1)\\
\eta (x,t_1)&=h(x)\quad\text{ in }\Omega
\endaligned\right.\tag 2.3
$$  
where
$$
H(x,t)=\left\{\aligned
&\lambda f(x)\left(\frac{g(u_1)^{-1}-g(u_2)^{-1}}{u_1-u_2}\right)\quad
\text{ if }u_1(x,t)\ne u_2(x,t)\\
&\lambda f(x)\left(\frac{1}{g}\right)'(u_1)\qquad\qquad\qquad
\text{ if }u_1(x,t)=u_2(x,t).\endaligned\right.\tag 2.4
$$
Then
$$\align
&\int_{\Omega}(u_1-u_2)(x,t_1)h(x)\,dx-\int_{\Omega}(u_{0,1}-u_{0,2})\eta\,dx\\
=&\int_0^{t_1}\int_{\Omega}\frac{\1}{\1 t}[(u_1-u_2)\eta]\,dx\,dt\\
=&\int_0^{t_1}\int_{\Omega}[(u_1-u_2)_t\eta+(u_1-u_2)\eta_t]\,dx\,dt\\
\le&\int_0^{t_1}\int_{\Omega}[\eta\Delta (u_1-u_2)
+\lambda\eta(g(u_1)^{-1}-g(u_2)^{-1})f +(u_1-u_2)\eta_t]\,dx\,dt\\
=&\int_0^{t_1}\int_{\Omega}(u_1-u_2)[\eta_t+\Delta\eta+H\eta]\,dx\,dt\\
=&0.
\endalign
$$
Hence
$$
\int_{\Omega}(u_1-u_2)(x,t_1)h(x)\,dx
\le\int_{\Omega}(u_{0,1}-u_{0,2})\eta\,dx.\tag 2.5
$$
Let $b=\sup_{\Omega\times (0,T)}|H(x,t)|$. By (2.1), (2.4) and either (1.10)
or (2.2), $b<\infty$. By the maximum principle $\eta\ge 0$. By (2.3),
$$\align
&\eta_t+\Delta\eta+b\eta\ge 0\qquad\quad\text{in }\Omega\times (0,t_1)\\
&(e^{bt}\eta)_t+\Delta (e^{bt}\eta)\ge 0\quad\text{ in }\Omega\times (0,t_1).
\endalign
$$
Hence by the maximum principle,
$$
\eta (x,0)\le\max_{\2{\Omega}} (e^{bt_1}\eta (x,t_1))=e^{bt_1}
\|h\|_{L^{\infty}}\le e^{bt_1}.\tag 2.6
$$
By (2.5) and (2.6),
$$
\int_{\Omega}(u_1-u_2)(x,t_1)h(x)\,dx
\le e^{bt_1}\int_{\Omega}(u_{0,1}-u_{0,2})_+\,dx.\tag 2.7
$$
Let $A=\{x\in\Omega:u_1(x,t_1)>u_2(x,t_1)\}$. We now choose a sequence of 
function $h_k\in C_0^{\infty}(\Omega)$, $0\le h_k\le 1$, such that
$h_k\to\chi_A$ a.e. as $k\to\infty$. Putting $h=h_k$ in (2.7) and letting
$k\to\infty$,
$$
\int_{\Omega}(u_1-u_2)(x,t_1)_+\,dx
\le e^{bt_1}\int_{\Omega}(u_{0,1}-u_{0,2})_+\,dx.
$$
Since $t_1\in (0,T)$ is arbitrary, (i) follows. Similarly if both
$u_1$ and $u_2$ are solutions of (0.3) in $\Omega\times (0,T)$
with initial value $u_0=u_{0,1}, u_{0,2}$, respectively, then
$$
\int_{\Omega}(u_1-u_2)(x,t)_-\,dx
\le e^{bt_1}\int_{\Omega}(u_{0,1}-u_{0,2})_-\,dx\quad\forall 0<t<T.
\tag 2.8
$$
By (i) and (2.8), (ii) follows. 
\enddemo

\proclaim{Corollary 2.2}
Let $u_{0,1}, u_{0,2}\in L^1(\Omega)$ be such that $u_{0,1}\le u_{0,2}$
in $\Omega$. Let $f\in C(\2{\Omega})$ and $0<g\in C^2((-\infty,1))$.
Suppose $u_1$, $u_2$, are the subsolution and supersolution 
of (0.3) in $\Omega\times (0,T)$ with initial value $u_0=u_{0,1}, u_{0,2}$, 
respectively. Suppose (2.1) holds and either (1.10) holds or (2.2) holds
for some constant $a_2<1$. Then $u_1\le u_2$ in $\2{\Omega}\times (0,T)$.
\endproclaim

\proclaim{Corollary 2.3}
Let $u_0\in L^1(\Omega)$, $f\in C(\2{\Omega})$ and $0<g\in C^2((-\infty,1))$
satisfy (1.10). Then the solution of (0.3) in $\Omega\times (0,T)$  is unique.
\endproclaim

\proclaim{Corollary 2.4}
Let $u_0\in L^1(\Omega)$, $f\in C(\2{\Omega})$ and $0<g\in C^2((-\infty,1))$. 
Then the solution of (0.3) in $\Omega\times (0,T)$ is unique in the class 
of functions on $\2{\Omega}\times (0,T)$ which are uniformly bounded below 
on $\2{\Omega}\times (0,T']$ for any $0<T'<T$.
\endproclaim

\proclaim{Theorem 2.5}
Let $u_0$ satisfy (0.4) for some constant $0<a<1$. Let $f$ satisfy
(0.1) and $g$ satisfy (0.2). Then for any $\lambda\ge 0$ there exists 
$T>0$ such that (0.3) has a solution which satisfies
$$
u(x,t)=\int_{\Omega}G(x,y,t)u_0(y)\,dy+\lambda\int_0^t\int_{\Omega}
G(x,y,t-s)\frac{f(y)}{g(u(y,s))}\,dy\,ds\quad\forall x\in\2{\Omega},
0<t<T,\tag 2.9
$$
$\Omega\times (0,T)$.
\endproclaim
\demo{Proof}
When $\lambda=0$, (0.3) reduces to the heat equation and the theorem follows 
from standard theory for heat equation \cite{F}. We next assume that
$\lambda>0$.
We divide the proof into two cases.

\noindent $\underline{\text{\bf Case 1}}$: $u_0\in C_0^{\infty}(\Omega)$
and $u_0$ satisfies (0.4) for some constant $0<a<1$.

Let 
$$
T=\frac{(1-a)}{4\lambda\|f\|_{L^{\infty}}}g((1+a)/2),\tag 2.10
$$
$$
w(x,t)=\int_{\Omega}G(x,y,t)u_0(y)\,dy,\tag 2.11
$$
and
$$
u_1(x,t)=w(x,t)+\lambda\int_0^t\int_{\Omega}
G(x,y,t-s)\frac{f(y)}{g(u_0(y))}\,dy\,ds\quad\forall x\in\2{\Omega},
0<t<T.\tag 2.12
$$ 
Then $w$ satisfies
$$\left\{\aligned 
\1_tw&=\Delta w\quad\text{ in }\Omega\times (0,\infty)\\
w(x,t)&=0\qquad\forall x\in\1\Omega, t>0\\
w(x,0)&=u_0(x)\quad\text{ in }\Omega.\endaligned\right.\tag 2.13
$$
Let $T_1=\sup\{0<t_1<T:u_1(x,t)<(1+a)/2\quad\forall x\in\2{\Omega}, 
0<t\le t_1\}$. Suppose $T_1<T$. By (0.2), (0.4), (0.6), (2.10) and (2.12), 
$\forall x\in\2{\Omega}, 0<t\le T_1$,
$$
u_1(x,t)\le a+\lambda (\|f\|_{L^{\infty}}/g(a))t
\le a+\frac{(1-a)}{4}\frac{g((1+a)/2)}{g(a)}<\frac{1+a}{2}.
$$
By continuity of $u_1$ there exists $0<\delta<(T-T_1)/2$ such that
$$
u_1(x,t)<\frac{1+a}{2}
$$
holds for all $x\in\2{\Omega}$, $0<t\le T_1+\delta$. This contradicts
the maximality of $T_1$. Hence $T_1=T$ and (2.12) holds for all 
$x\in\2{\Omega}$, $0<t\le T$. Suppose $u_1,u_2,\dots,u_k$, are defined. 
We define
$$
u_{k+1}(x,t)=w(x,t)+\lambda\int_0^t\int_{\Omega}G(x,y,t-s)
\frac{f(y)}{g(u_k(y,s))}\,dy\,ds\quad\forall x\in\2{\Omega},0<t<T.\tag 2.14
$$
Let $T_k=\sup\{0<t_1<T:u_k(x,t)<(1+a)/2\quad\forall x\in\2{\Omega}, 
0<t\le t_1\}$. We claim that $T_k=T$ for all $k\in\Bbb{Z}^+$.
We will prove this claim by induction. Note that $T_1=T$ is already proved
before. Suppose $T_1=T_2=\cdots=T_k=T$ but $T_{k+1}<T$. Then
$$
u_k(x,t)<\frac{1+a}{2}\quad\forall x\in\2{\Omega}, 0<t<T.\tag 2.15
$$ 
By (0.2), (0.4), (2.10), (2.14) and (2.15), 
$$
u_{k+1}(x,t)\le a +\lambda (\|f\|_{L^{\infty}}/g((1+a)/2))t
\le a+\frac{(1-a)}{4}<\frac{1+a}{2}\quad\forall x\in\2{\Omega},0<t\le T_1.
$$
By continuity of $u_{k+1}$ there exists $0<\delta<(T-T_{k+1})/2$ such that
$$
u_{k+1}(x,t)<\frac{1+a}{2}\tag 2.16
$$
holds for all $x\in\2{\Omega}$, $0<t\le T_{k+1}+\delta$. This contradicts
the maximality of $T_{k+1}$. Hence $T_{k+1}=T$ and (2.16) holds for all 
$x\in\2{\Omega}$, $0<t\le T$. Thus by induction $T_k=T$ for all 
$k\in\Bbb{Z}^+$. Hence (2.15) holds for all $k\in\Bbb{Z}^+$.
Since 
$$
w(\cdot,t)\to u_0\quad\text{ in }L^1(\Omega)\quad\text{ as }t\to 0,
\tag 2.17
$$ 
by (0.6), (2.14) and (2.15),
$$
u_k(\cdot,t)\to u_0\quad\text{ in }L^1(\Omega)\quad\text{ as }t\to 0.
\tag 2.18
$$
By (2.12) and (2.14),
$$
w(x,t)\le u_k(x,t)\quad\text{ in }\Omega\times (0,T)\quad\forall 
k\in\Bbb{Z}^+.\tag 2.19
$$
By (2.12), $u_1$ is continuously differentiable in $x$ and $t$. Then by
(2.14), (2.19) and standard parabolic theory \cite{F}, $u_k\in C^{2,1}
(\2{\Omega}\times (0,T])$ for all $k\ge 2$. Then by (2.14), (2.15) and 
(2.18), (2.19), $\forall k\ge 2$, $u_k$ satisfies
$$\left\{\aligned
\frac{\1 u_k}{\1 t}-\Delta u_k&
=\frac{\lambda f}{g(u_{k-1})}\quad\text{ in }
\Omega\times (0,T)\\
u_k(x,t)&=0\qquad\quad\text{ on }\1\Omega\times (0,T)\\
u_k(x,0)&=u_0(x)\quad\text{ in }\Omega.
\endaligned\right.\tag 2.20
$$
By (2.15), (2.19), (2.20) and the parabolic Schauder estimates \cite{LSU}, 
the sequence $\{u_k\}_{k=2}^{\infty}$ are uniformly Holder continuous 
on $\2{\Omega}\times [0,T]$. Then by (2.15), (2.19), (2.20) and the 
Schauder estimates for the heat equation (\cite{F},\cite{LSU}) 
$\{u_k\}_{k=2}^{\infty}$ are uniformly bounded in
$C^{2+\beta,1+(\beta/2)}(K)$ for any compact subset 
$K\subset\2{\Omega}\times (0,T]$ where $0<\beta<1$ is some constant. 
By the Ascoli theorem and a 
diagonalization argument $\{u_k\}_{k=2}^{\infty}$ has a subsequence
which we may assume without loss of generality to be the sequence itself
which converges uniformly in $C^{2+\beta,1+(\beta/2)}(K)$ to some 
function $u$ for any
compact subset $K\subset\2{\Omega}\times (0,T]$ as $k\to\infty$. Then by 
(2.14), (2.15), (2.19) and (2.20) $u$ satisfies (2.9),
$$
w(x,t)\le u(x,t)\le\frac{1+a}{2}\quad\forall x\in\2{\Omega}, 0<t\le T,
\tag 2.21
$$
and
$$\left\{\aligned
\frac{\1 u}{\1 t}-\Delta u&
=\frac{\lambda f}{g(u)}\quad\text{ in }
\Omega\times (0,T)\\
u(x,t)&=0\qquad\quad\text{on }\1\Omega\times (0,T)\\
\endaligned\right.\tag 2.22
$$  
By (0.6), (2.9), (2.17) and (2.21), $u$ satisfies (0.5). Hence $u$ is a 
solution of (0.3) in $\Omega\times (0,T)$.

\noindent $\underline{\text{\bf Case 2}}$: $u_0$ satisfies (0.4) for some
constant $0<a<1$.

We choose a sequence of function $\{u_{0,k}\}_{k=1}^{\infty}\in 
C_0^{\infty}(\Omega)$ such that $u_{0,k}$ converges to $u_0$ in $L^1(\Omega)$
and a.e. as $k\to\infty$. For any $k\in\Bbb{Z}^+$, by case 1 there exists
a solution $u_k$ of (0.3) in $\Omega\times (0,T)$ with initial value $u_{0,k}$
which satifies
$$
u_k(x,t)=\int_{\Omega}G(x,y,t)u_{0,k}(y)\,dy
+\lambda\int_0^t\int_{\Omega}G(x,y,t-s)
\frac{f(y)}{g(u_k(y,s))}\,dy\,ds\tag 2.23
$$
for any $x\in\2{\Omega},0<t<T$, and
$$
w(x,t)\le u_k(x,t)\le\frac{1+a}{2}\quad\forall x\in\2{\Omega}, 0<t\le T.
\tag 2.24
$$
Since $\{u_k\}_{k=1}^{\infty}$ satisfy (0.3)  with initial value $u_{0,k}$
in $\Omega\times (0,T)$, by the parabolic Schauder estimates \cite{LSU}, 
the sequence $\{u_k\}_{k=1}^{\infty}$ 
are uniformly Holder continuous on $\2{\Omega}\times (\delta_1,T]$ for any
$0<\delta_1<T$. Then by the parabolic Schauder estimates 
(\cite{F},\cite{LSU}) $\{u_k\}_{k=1}^{\infty}$ are uniformly bounded in
$C^{2+\beta,1+(\beta/2)}(K)$ for any compact subset 
$K\subset\2{\Omega}\times (0,T]$
where $0<\beta<1$ is some constant. By the Ascoli theorem and a 
diagonalization argument $\{u_k\}_{k=1}^{\infty}$ has a subsequence
which we may assume without loss of generality to be the sequence itself
which converges uniformly in $C^{2+\beta,1+(\beta/2)}(K)$ 
to some function $u$ for any
compact subset $K\subset\2{\Omega}\times (0,T]$ as $k\to\infty$. Then $u$ 
satisfies (2.22). Letting $k\to\infty$ in (2.23) and (2.24), we get (2.9) 
and (2.21). By (0.6), (2.9), (2.17) and (2.21), $u$ satisfies (0.5). 
Hence $u$ is 
a solution of (0.3) in $\Omega\times (0,T)$ and the theorem follows.
\enddemo

By Corollary 2.2 and the Duhamel principle we have the following corollary.

\proclaim{Corollary 2.6}
Let $f$ satisfy (0.1), $g$ satisfy (0.2) and $u_0$ satisfy (0.4) for 
some constant $a<1$. Suppose $u$ is a bounded solution of (0.3) in 
$\Omega\times (0,T)$. Then $u$ satisfies (2.9) in $\2{\Omega}\times (0,T)$.
\endproclaim

\proclaim{Corollary 2.7}
Let $f$ satisfy (0.1) and $g$ satisfy (0.2).
Let $u_{0,1}, u_{0,2}\in L^{\infty}(\Omega)$ be such that $u_{0,1}
\le u_{0,2}\le a<1$ for some constant $0<a<1$ and $u_{0,1}\not\equiv u_{0,2}$. 
Suppose $u_1$, $u_2$, are bounded solutions of (0.3) in $\Omega\times (0,T)$ 
with initial values $u_{0,1}, u_{0,2}$ respectively. Then
$$
u_1<u_2\quad\text{ in }\Omega\times (0,T).
$$ 
\endproclaim
\demo{Proof}
By Corollary 2.2 $u_1\le u_2$ in $\Omega\times (0,T)$. By Corollary 2.6 both
$u_1$ and $u_2$ satisfies (2.9) with $u_0=u_{0,1}, u_{0,2}$ respectively.
By (2.9) for $u_1$, $u_2$, (0.2) and the positivity of the Green function
for the heat equation the corollary follows. 
\enddemo

By Corollary 2.2, Theorem 2.5 and a continuity argument we have the 
following theorem.

\proclaim{Theorem 2.8}
Let $f$ satisfy (0.1) and $g$ satisfy (0.2).
Let $\lambda\ge 0$. Suppose ($S_{\lambda}$) has a supersolution
$v_{\lambda}$. Let $u_0\in L^{\infty}(\Omega)$ satisfy 
$$
u_0\le v_{\lambda}\quad\text{ in }\Omega.
$$
Then (0.3) has a unique bounded global solution which satisfies (2.9) and
$$
\inf_{\Omega}u_0\le u(x,t)\le v_{\lambda}(x)\quad\forall\2{\Omega}
\times (0,\infty).
$$
\endproclaim

\proclaim{Theorem 2.9}
Let $f$ satisfy (0.1) and $g$ satisfy (0.2).
Let $0\le\lambda\le\lambda^{\ast}$ and let $v_{\lambda}$ be a supersolution
of ($S_{\lambda}$). Let $u_0\in L^1(\Omega)$ satisfy 
$$
u_0\le v_{\lambda}\quad\text{ in }\Omega.
$$
Then (0.3) has a global solution $u$ which satisfies (2.9) and
$$
w(x,t)\le u(x,t)\le v_{\lambda}(x)\quad\forall\2{\Omega}\times (0,\infty)
\tag 2.25
$$
where $w$ is given by (2.11). The solution is unique within the family of
functions satisfying (2.25) if either (1.10) holds
or
$$
\sup_{s\le a}\left(\frac{1}{g}\right)'(s)<\infty\quad\forall a<1.\tag 2.26
$$
\endproclaim
\demo{Proof}
For any $k\in\Bbb{Z}^+$, let $u_{0,k}=\max (u_0,-k)$. Then
$$
u_{0,k+1}\le u_{0,k}\quad\text{ and }-k\le u_{0,k}\le v_{\lambda}
\quad\text{ in }\Omega\quad\forall k\in\Bbb{Z}^+.
$$
By Corollary 2.2 and Theorem 2.8 for any $k\in\Bbb{Z}^+$ there exists a 
global bounded solution $u_k$ of (0.3) with initial value $u_{0,k}$ which 
satisfies (2.23) in $\Omega\times (0,\infty)$,
$$
-k\le u_k\le v_{\lambda}\quad\text{ in }\Omega\times (0,\infty)
\quad\forall k\in\Bbb{Z}^+.\tag 2.27
$$
and
$$
u_{k+1}\le u_k\quad\text{ in }\Omega\times (0,\infty)\quad\forall 
k\in\Bbb{Z}^+.\tag 2.28
$$
By (2.23),
$$
w_k(x,t)\le u_k\quad\text{ in }\Omega\times (0,\infty)\quad\forall 
k\in\Bbb{Z}^+\tag 2.29
$$
where
$$
w_k(x,t)=\int_{\Omega}G(x,y,t)u_{0,k}(y)\,dy\tag 2.30
$$
is the solution of (2.13) with initial value $u_{0,k}$. Let $w$ be given by 
(2.11). Since $|u_{0,k}|\le |u_0|$ in $\Omega$, by (0.6), (2.30) and 
the Lebesgue 
dominated convergence theorem $w_k$ converges uniformly to $w$ on 
$\2{\Omega}\times[\delta_1,\infty)$ as $k\to\infty$ for any $\delta_1>0$. 
Hence by (2.27) and (2.29), the sequence $\{u_k\}_{k=1}^{\infty}$ are 
uniformly bounded on $\2{\Omega}\times[\delta_1,\infty)$ for any $\delta_1>0$.
Since $u_k$ satisfies (0.3) in $\Omega\times (0,\infty)$ with initial value
$u_{0,k}$, by the Schauder estimates \cite{LSU} 
$\{u_k\}_{k=1}^{\infty}$ are uniformly bounded in 
$C^{2+\beta,1+(\beta/2)}(\2{\Omega}\times[\delta_1,\infty))$ 
for any $\delta_1>0$ where $0<\beta<1$ is some constant. By (2.28), 
the Ascoli theorem and a diagonalization argument 
$\{u_k\}_{k=1}^{\infty}$ has a subsequence which we may assume without 
loss of generality to be the sequence itself which decreases and converges 
uniformly in $C^{2+\beta,1+(\beta/2)}(\2{\Omega}\times[\delta_1,\infty))$ 
to some function $u$ for any $\delta_1>0$ as $k\to\infty$.
 
Then $u$ satisfies (2.22) and (2.25). Letting $k\to\infty$ in (2.23) we get 
(2.9). By (2.9) and (2.17)  $u$ satisfies (0.5). Hence $u$ is 
a solution of (0.3) in $\Omega\times (0,T)$. If (1.10) holds, by Corollary 
2.3 the solution is unique. 

Suppose (2.26) holds. Suppose $u_1$, $u_2$, are both solutions of (0.3)
in $\Omega\times (0,\infty)$. Then by (2.25) and the Duhamel principle, 
both $u_1$, $u_2$, satisfies (2.9). Putting $u=u_1, u_2$, in (2.9) and 
subtracting the resulting equations, we get
$$\align
u_1(x,t)-u_2(x,t)&=\lambda\int_0^t\int_{\Omega}G(x,y,t-s)f(y)
\left(\frac{1}{g(u_1(y,s))}-\frac{1}{g(u_2(y,s))}\right)\,dy\,ds\\
&\le\lambda\|f\|_{L^{\infty}}\int_0^t\int_{\Omega}G(x,y,t-s)
\left(\frac{1}{g}\right)'(\xi (y,s))(u_1(y,s)-u_2(y,s))_+\,dy\,ds\\
&\le a_0\lambda\|f\|_{L^{\infty}}\int_0^t\int_{\Omega}G(x,y,t-s)
(u_1(y,s)-u_2(y,s))_+\,dy\,ds\\
&\le a_0\lambda\|f\|_{L^{\infty}}t\sup_{\Omega\times (0,T)}
(u_1-u_2)_+\quad\forall x\in\Omega, 0<t<T
\endalign
$$ 
for any $T>0$ where $\xi(y,s)$ is some number between $u_1(y,s)$ and 
$u_2(y,s)$,
$$
a_0=\sup_{s\le\|v_{\lambda}\|_{L^{\infty}}}\left(\frac{1}{g}\right)'(s).
$$ 
Hence
$$
\sup_{\Omega\times (0,T)}(u_1-u_2)_+
\le a_0\lambda\|f\|_{L^{\infty}}T\sup_{\Omega\times (0,T)}
(u_1-u_2)_+.\tag 2.31
$$ 
We now choose $T=1/(1+2a_0\lambda\|f\|_{L^{\infty}})$. Then by (2.31),
$$
\sup_{\Omega\times (0,T)}(u_1-u_2)_+=0\quad
\Rightarrow\quad u_1\le u_2\quad\text{ in }\2{\Omega}\times (0,T).
$$
By interchanging the role of $u_1$ and $u_2$ we get
$$
u_2\le u_1\quad\text{ in }\2{\Omega}\times (0,T).
$$
Hence
$$
u_1=u_2\quad\text{ in }\2{\Omega}\times (0,T).
$$
By dividing the time interval into disjoint intervals of length $T$ and 
repeating the above argument we get
$$
u_1=u_2\quad\text{ in }\2{\Omega}\times (0,\infty)
$$
and the theorem follows. 
\enddemo

\proclaim{Theorem 2.10}
Let $g$ satisfy (0.2) and 
$$
0\le f\in C^{\alpha}(\Bbb{R}^n)\quad\text{ for some constant }0<\alpha<1.
$$ 
Let $u_0\in L^1(\Bbb{R}^n)$ be such that $u_0\le a$ in $\Bbb{R}^n$ for
some constant $a<1$. Then for any $\lambda\ge 0$ there exists a constant 
$T>0$ such that the Cauchy problem
$$\left\{\aligned
&u_t=\Delta u+\frac{\lambda f(x)}{g(u)}\quad\text{ in }\Bbb{R}^n\times (0,T)\\
&u(x,0)=u_0\qquad\quad\text{ in }\Bbb{R}^n
\endaligned\right.\tag 2.32
$$
has a solution $u$ which satisfies
$$
u(x,t)=\int_{\Bbb{R}^n}Z(x,y,t)u_0(y)\,dy+\lambda\int_0^t\int_{\Bbb{R}^n}
Z(x,y,t-s)\frac{f(y)}{g(u(y,s))}\,dy\,ds
\tag 2.33
$$
in $\Bbb{R}^n\times (0,T)$ where $Z(x,y,t)=(4\pi)^{-\frac{n}{2}}
e^{-|x-y|^2/4t}$.
\endproclaim
\demo{Proof}
If $\lambda=0$ or $f\equiv 0$ in $\Bbb{R}^n$, (2.32) reduces to the heat 
equation and the result follows by standard results on heat equation 
\cite{F}. Hence we may assume without loss of generality that $\lambda>0$
and $f\not\equiv 0$ in $\Bbb{R}^n$. Let $T$ be given by (2.10). For any 
$R>0$ let $G_R(x,y,t)$ be the Dirichlet Green function of the heat equation 
in $B_R\times (0,\infty)$. By the proof of Theorem 2.5 for any $k\ge 1$ 
there exists a solution $u_k$ of
$$\left\{\aligned
&u_t=\Delta u+\frac{\lambda f(x)}{g(u)}\quad\text{ in }B_k\times (0,T)\\
&u(x,t)=0\qquad\qquad\text{ on }\1 B_k\times (0,T)\\
&u(x,0)=u_0\qquad\quad\text{ in }B_k
\endaligned\right.\tag 2.34
$$
which satisfies
$$
u_k(x,t)=\int_{B_k}G_k(x,y,t)u_0(y)\,dy+\lambda\int_0^t\int_{B_k}
G_k(x,y,t-s)\frac{f(y)}{g(u_k(y,s))}\,dy\,ds
\tag 2.35
$$
for any $(x,t)\in B_k\times (0,T)$ and
$$
w_k(x,t)\le u_k(x,t)\le\frac{1+a}{2}\quad\text{ in }B_k\times (0,T)
\quad\forall k\ge 1\tag 2.36
$$
where
$$
w_k(x,t)=\int_{B_k}G_k(x,y,t)u_0(y)\,dy
$$
Since $G_k(x,y,t)\le G_{k+1}(x,y,t)$ in $B_k\times (0,T)$ for any $k\ge 1$,
by the construction of solutions in Theorem 2.5,
$$
u_k\le u_{k+1}\quad\text{ in }B_k\times (0,T)\quad\forall k\ge 1.\tag 2.37
$$
Since $w_k$ converges uniformly to 
$$
w(x,t)=\int_{\Bbb{R}^n}Z(x,y,t)u_0(y)\,dy\tag 2.38
$$
as $k\to\infty$, by (2.36) the sequence $\{u_k\}_{k=1}^{\infty}$ is 
uniformly bounded on every compact subset of $\Bbb{R}^n\times (0,T)$. 
By (2.34) for $u_k$, (2.36), and the parabolic Schauder estimates the 
sequence $\{u_k\}_{k=1}^{\infty}$ is uniformly Holder continuous on every 
compact subset of $\Bbb{R}^n\times (0,T)$. Then by (2.34) for $u_k$, (2.36), 
and the parabolic Schauder estimates the sequence 
$\{u_k\}_{k=1}^{\infty}$ is uniformly bounded in 
$C^{2+\beta,1+(\beta/2)}(K)$ for any compact subset $K\subset\Bbb{R}^n\times 
(0,T)$ where $0<\beta<1$ is some constant. Then by (2.35), (2.36), (2.37), the 
Ascoli Theorem and a diagonalization argument the sequence 
$\{u_k\}_{k=1}^{\infty}$ has a subsequence which we may
assume without loss of generality to be the sequence $\{u_k\}_{k=1}^{\infty}$
itself which increases and converges uniformly in 
$C^{2+\beta,1+(\beta/2)}(K)$ for any compact subset $K\subset\Bbb{R}^n\times 
(0,T)$ to a function $u$ which satisfies (2.33),
$$
u_t=\Delta u+\frac{\lambda f(x)}{g(u)}\quad\text{ in }\Bbb{R}^n\times (0,T),
$$
and
$$
w(x,t)\le u(x,t)\le\frac{1+a}{2}\quad\text{ in }\Bbb{R}^n\times (0,T)
\tag 2.39
$$
Since $w(x,t)\to u_0$ as $t\to 0$, by (2.33) and (2.39) $u(x,t)\to u_0$ as 
$t\to 0$. Hence $u$ satisfies (2.32) in $\Bbb{R}^n\times (0,T)$.
\enddemo

$$
\text{\bf Section 3}
$$

In this section we will prove the convergence of solutions of (0.3) for any
$0\le\lambda<\lambda^{\ast}$ as $t\to\infty$. We also obtain various conditions
for the solutions of (0.3) to have finite touchdown time.

\proclaim{Theorem 3.1}
Suppose $f$ satisfies (0.1) and $g$ satisfies (0.2). Let $0\le\lambda
<\lambda^{\ast}$ and let $v_{\lambda}$ be the unique minimal solution 
of ($S_{\lambda}$) given by Theorem 1.1. Let $u_0$ satisfies 
$$
u_0\le v_{\lambda}\quad\text{ in }\Omega
$$
and let $u$ be the global solution of (0.3) constructed in Theorem 2.9.
Then $u$ converges uniformly on $\2{\Omega}$ to $v_{\lambda}$ as $t\to\infty$.
\endproclaim
\demo{Proof}
Note that the theorem is proved by N.~Ghoussoub and Y.~Guo in \cite{GhG2} 
for the case $g(s)=(1-s)^2$ and $u_0=0$ in $\Omega$ and by T.~Suzuki, etc. in 
\cite{KMS} for the case $g(s)=(1-s)^p$ and 
$0\le u_0\le v_{\lambda}$ in $\Omega$. Both are based on proving the 
positivity of $u_t$ in $\Omega\times (0,\infty)$ when $u_0=0$ using a 
modification of Fujita's technique \cite{Fu}. This approach is not 
applicable in our case and we will use a different proof for the 
convergence result. 

By Theorem 2.9 $u$ satisfies (2.9) and (2.25) with $w$ being given by (2.11).
Let $\{t_k\}_{k=1}^{\infty}$, $t_k\ge 1$ for all $k\ge 1$, be a sequence 
such that $t_k\to\infty$ as $k\to\infty$. By (2.25) and the parabolic 
Schauder estimates \cite{LSU} $u(x,t)$ is uniformly bounded in 
$C^{2+\beta,1+(\beta/2)}(\2{\Omega}\times [1,\infty))$ where $0<\beta<1$ 
is some constant. Then by the Ascoli theorem $\{t_k\}_{k=1}^{\infty}$ has a 
subsequence $\{t_{i_k}\}_{k=1}^{\infty}$ such that 
$u(x,t_{i_k}+t)$ converges uniformly in $C^{2,1}(\2{\Omega}\times [0,1])$ 
to some function $v_1$ as $k\to\infty$. Let $v(x)=v_1(x,0)$. Multiplying 
(0.3) by $u_t$ and integrating over $\Omega\times (1,t)$,
$$\align
\int_1^t\int_{\Omega}u_t^2\,dx\,dt
=&\int_1^t\int_{\Omega}u_t\Delta u\,dx\,dt
+\lambda\int_1^t\int_{\Omega}\frac{f(x)u_t}{g(u)}\,dx\,dt\\
=&-\frac{1}{2}\int_1^t\frac{\1}{\1 t}\left(\int_{\Omega}|\nabla u|^2\,dx
\right)\,dt
+\lambda\int_{\Omega}f(x)\left(\int_{u(x,1)}^{u(x,t)}\frac{ds}{g(s)}
\right)\,dx\\
\le&\frac{1}{2}\int_{\Omega}|\nabla u(x,1)|^2\,dx
+\lambda\|f\|_{L^{\infty}}|\Omega|(a_2-a_1)\max_{a_1\le s\le a_2}(1/g(s))
\endalign
$$
holds for all $t\ge 1$ where $a_1=\min_{\2{\Omega}}u(x,1)$, 
$a_2=\max_{\2{\Omega}}v_{\lambda}$. Letting $t\to\infty$,
$$
\int_1^{\infty}\int_{\Omega}u_t^2\,dx\,dt
\le\frac{1}{2}\int_{\Omega}|\nabla u(x,1)|^2\,dx
+\lambda\|f\|_{L^{\infty}}|\Omega|(a_2-a_1)\max_{a_1\le s\le a_2}(1/g(s)).
$$
Hence
$$
\int_{t_{i_k}}^{t_{i_k}+1}\int_{\Omega}u_t^2\,dx\,dt\to 0\quad\text{ as }
k\to\infty.
$$
Thus
$$\align
\int_{\Omega}|u(x,t_{i_k}+t)-u(x,t_{i_k})|\,dx
\le&\int_{t_{i_k}}^{t_{i_k}+1}\int_{\Omega}|u_t|\,dx\,dt\\
\le&|\Omega|^{\frac{1}{2}}
\left(\int_{t_{i_k}}^{t_{i_k}+1}\int_{\Omega}u_t^2\,dx\,dt
\right)^{\frac{1}{2}}\\
\to&0\qquad\qquad\text{ as }k\to\infty\\
\Rightarrow\qquad\qquad\int_{\Omega}|v_1(x,t)-v(x)|\,dx
=&0\quad\forall 0\le t\le 1\\
\Rightarrow\qquad\qquad\qquad\qquad\qquad 
v_1(x,t)=&v(x)\quad\forall x\in\2{\Omega},0\le t\le 1.
\endalign
$$
Hence $u(x,t_{i_k}+t)$ converges uniformly to $v(x)$ on $\2{\Omega}\times
[0,1]$ as $k\to\infty$. Putting $t=t_{i_k}$ and letting $k\to\infty$ in (2.25),
$$
0\le v(x)\le v_{\lambda}(x)\quad\text{ in }\2{\Omega}.\tag 3.1
$$
Integrating (0.3) over $(t_{i_k}, t_{i_k}+1)$,
$$
u(x,t_{i_k}+1)-u(x,t_{i_k})=\int_{t_{i_k}}^{t_{i_k}+1}
\Delta u(x,s)\,ds+\int_{t_{i_k}}^{t_{i_k}+1}
\frac{\lambda f(x)}{g(u(x,s))}\,ds\quad\text{ on }\2{\Omega}.
$$
Letting $k\to\infty$ we get that $v$ satisfies ($S_{\lambda}$). 
Since $v_{\lambda}$ is the minimal solution of ($S_{\lambda}$), by (3.1),
$$
v(x)=v_{\lambda}(x)\quad\text{ on }\2{\Omega}.
$$
Since the sequence $\{t_k\}_{k=1}^{\infty}$ is arbitrary, $u(x,t)$ 
converges uniformly to $v_{\lambda}$ on $\2{\Omega}$ as $t\to\infty$ and 
the theorem follows.
\enddemo

By (ii) of Theorem 1.1 and an argument similar to the proof of 
Theorem 3.1 we have the following theorem.

\proclaim{Theorem 3.2}
Suppose $f$ satisfies (0.1) and $g$ satisfies (0.2). Let $\lambda
>\lambda^{\ast}$ and let $u$ be a solution of (0.3). Then 
either $T_{\lambda}<\infty$ or $u$ touchdowns at time infinity.
\endproclaim

\proclaim{Theorem 3.3}
Let $f$ satisfy (0.1) and (1.4) and $g$ satisfy (0.2) and (1.10). 
Let $\lambda_1=(\mu_1/\delta_1)\sup_{0\le s\le 1}sg(s)$. Then for 
any solution $u$ of (0.3) with initial value $u_0$ and $\lambda>\lambda_1$, 
we have 
$$
T_{\lambda}\le\frac{1}{(\lambda-\lambda_1)\delta_1}
\int_{E(0)}^1g(s)\,ds\tag 3.2
$$
where $E(0)=\int_{\Omega}u_0\phi_1\,dx$.
Moreover if $g$ also satisfies
$$
g(s)\to 0\quad\text{ as }s\nearrow 1,\tag 3.3
$$
then there exists a constant $a_0<1$ such that if
$$
\int_{\Omega}u_0\phi_1\,dx\ge a_0,\tag 3.4
$$
then for any solution $u$ of (0.3) with $\lambda>0$ and initial 
value $u_0$ we have $T_{\lambda}\le (1-a_0)/10$.
\endproclaim
\demo{Proof}
We will use a modification of the argument of \cite{GPW} and \cite{KMS}
to prove the theorem. Suppose $u$ is a solution of (0.3) with $\lambda>0$
and initial value $u_0$. Let 
$$
E(t)=\int_{\Omega}u(x,t)\phi_1(x)\,dx.
$$
Multiplying (0.3) by $\phi_1$ and integrating over $\Omega$, by the Green
theorem, (1.10) and the Jensen inequality,
$$\align
\frac{d}{dt}E(t)=&\frac{d}{dt}\left(\int_{\Omega}u\phi_1\,dx\right)
=\int_{\Omega}\phi_1\Delta u\,dx+\lambda\int_{\Omega}\frac{f\phi_1}{g(u)}\,dx\\
\ge&-\mu\int_{\Omega}u\phi_1\,dx+\lambda\delta_1
\int_{\Omega}\frac{\phi_1}{g(u)}\,dx\\
\ge&-\mu E(t)+\lambda\frac{\delta_1}{g(E(t))}\tag 3.5
\endalign
$$  
Note $E(t)\le 1$ for any $t>0$. We now divide the proof into two cases.

\noindent $\underline{\text{\bf Case 1}}$: $\lambda>\lambda_1$. 

Then the right hand side is 
$$
\ge (\lambda-\lambda_1)\frac{\delta_1}{g(E(t))}.\tag 3.6
$$
Integrating (3.5) over $(0,t)$, by (3.6),
$$
t\le\frac{1}{(\lambda-\lambda_1)\delta_1}\int_{E(0)}^1g(s)\,ds
$$
and (3.2) follows.

\noindent $\underline{\text{\bf Case 2}}$: $\lambda>0$ and (3.3), (3.4), 
hold for some constant $a_0$ to be determined later. 

By (3.3) there exists a constant $a_0<1$ such that
$$
-\mu y+\lambda\frac{\delta_1}{g(y)}\ge 10\quad\forall a_0\le y<1.\tag 3.7
$$
Integrating (3.5) over $(0,t)$, by (3.4) and (3.7),
$$
10t\le E(t)-E(0)\le 1-E(0)\quad\Rightarrow\quad T_{\lambda}
\le\frac{1-E(0)}{10}\le\frac{1-a_0}{10}.
$$
\enddemo

By Corollary 2.2, Theorem 2.9, Theorem 3.3 and a comparison argument we 
have the following corollary.

\proclaim{Corollary 3.4}
Let $f$ satisfy (0.1),
$$
\delta_R=\inf_{B_R(x_0)}f>0
$$
for some $B_R(x_0)\subset\Omega$, and let $g$ satisfy (0.2) and 
(1.10). Let $\mu_R$ be the first eigenvalue of $-\Delta$ in $B_R(x_0)$
and let $\phi_R$ be the first positive eigenfunction of $-\Delta$ in 
$B_R(x_0)$ normalized such that $\int_{B_R(x_0)}\phi_R\,dx=1$.
Let $\lambda_R=(\mu_R/\delta_R)\sup_{0\le s\le 1}sg(s)$. Then for 
any solution $u$ of (0.3) with initial value $u_0\ge 0$ and 
$\lambda>\lambda_R$, we have 
$$
T_{\lambda}\le\frac{1}{(\lambda-\lambda_1)\delta_R}
\int_{E_1(0)}^1g(s)\,ds
$$
where $E_1(0)=\int_{B_R(x_0)}u_0\phi_R\,dx$. Moreover if $g$ also satisfies
(3.3), then there exists a constant $a_1<1$ such that if
$u_0\ge 0$ and 
$$
\int_{B_R(x_0)}u_0\phi_R\,dx\ge a_1,
$$
then for any solution $u$ of (0.3) with $\lambda>0$ and initial 
value $u_0$ we have $T_{\lambda}\le (1-a_1)/10$.
\endproclaim

\proclaim{Theorem 3.5}
Let $f$ satisfy (0.1), $g$ satisfy (0.2) and $\lambda>\lambda^{\ast}$. 
Suppose $u_0$ satisfies (0.4) for some constant $a<1$ and
$$
\2{u}_0\le u_0\quad\text{ in }\Omega\tag 3.8
$$
for some subsolution $\2{u}_0\in C^2(\Omega)\cap C(\2{\Omega})$ of 
($S_{\lambda}$). If $u$ is the unique bounded solution of (0.3),
then $T_{\lambda}<\infty$. 
\endproclaim
\demo{Proof}
Suppose $u$ is a global bounded solution of (0.3). Let $\2{u}$  be the 
unique bounded solution of (0.3) 
with initial value $\2{u}_0$ given by Theorem 2.5 and Corollary 2.2.
Then by Theorem 2.5, Corollary 2.2 and a continuity argument $\2{u}$ 
can be extended to a global solution of (0.3) with initial value 
$\2{u}_0$ which satisfies
$$
\2{u}\le u\quad\text{ in }\Omega\times (0,\infty).\tag 3.9
$$
By an argument similar
to the proof on P.4--6 of \cite{KMS} but with $(1-u)^p$ there being replaced 
by $g(u)$ we get that there exists a time $T>0$ such that
$$
\lim_{t\nearrow T}\sup_{\Omega}\2{u}(x,t)=1.\tag 3.10
$$
By (3.9) and (3.10), $\sup_{\Omega}u(x,t)$ will converges to $1$ before the 
time $T$. Hence $T_{\lambda}<\infty$.
\enddemo

\proclaim{Theorem 3.6}
Let $f$ satisfy (0.1) and $g$ satisfy (0.2). Let
$$
\lambda
>\mu_1\frac{\int_0^1g(s)\,ds}{\int_{\Omega}f\phi_1\,dx}
$$
and let $u$ be a solution of (0.3) with initial value $u_0$. Then 
$$
T_{\lambda}\le\frac{1}{(\lambda-\lambda')}\frac{\int_{\Omega}H(u_0)\phi_1\,dx}
{\int_{\Omega}f\phi_1\,dx}\tag 3.11
$$
where
$$
\lambda'=\mu_1\frac{\int_0^1g(s)\,ds}{\int_{\Omega}f\phi_1\,dx}.
$$
\endproclaim
\demo{Proof}
Let $H(u)$ be given by (1.3). Then
$$\align
\frac{d}{dt}\left(\int_{\Omega}H(u)\phi_1\,dx\right)
=&-\int_{\Omega}\phi_1 g(u)u_t\,dx\\
=&-\int_{\Omega}\phi_1 g(u)\Delta u\,dx
-\lambda\int_{\Omega}f\phi_1\,dx\\
=&\int_{\Omega}\phi_1 g'(u)|\nabla u|^2\,dx
+\int_{\Omega}g(u)\nabla\phi_1\cdot\nabla u\,dx
-\lambda\int_{\Omega}f\phi_1\,dx\\
\le&-\int_{\Omega}\nabla\phi_1\cdot\nabla H(u)\,dx
-\lambda\int_{\Omega}f\phi_1\,dx\\
\le&\int_{\Omega}H(u)\Delta\phi_1\,dx
-\int_{\1\Omega}H(u)\frac{\1\phi_1}{\1\nu}\,d\sigma
-\lambda\int_{\Omega}f\phi_1\,dx\\
\le&-\mu_1\int_{\Omega}H(u)\phi_1\,dx
-H(0)\int_{\1\Omega}\frac{\1\phi_1}{\1\nu}\,d\sigma
-\lambda\int_{\Omega}f\phi_1\,dx\\
=&-\mu_1\int_{\Omega}H(u)\phi_1\,dx
-H(0)\int_{\Omega}\Delta\phi_1\,dx
-\lambda\int_{\Omega}f\phi_1\,dx\\
=&-\mu_1\int_{\Omega}H(u)\phi_1\,dx
+\mu H(0)\int_{\Omega}\phi_1\,dx
-\lambda\int_{\Omega}f\phi_1\,dx\\
=&-\mu_1\int_{\Omega}H(u)\phi_1\,dx
+\mu_1H(0)-\lambda\int_{\Omega}f\phi_1\,dx\\
\le&-(\lambda-\lambda')\int_{\Omega}f\phi_1\,dx
\endalign
$$
Integrating over $(0,t)$,
$$
\int_{\Omega}H(u(x,t))\phi_1(x)\,dx\le\int_{\Omega}H(u_0)\phi_1\,dx
-(\lambda-\lambda')t\int_{\Omega}f\phi_1\,dx.
$$
Since the left hand side is positive while the right hand is negative for any 
$$
t>\frac{1}{(\lambda-\lambda')}\frac{\int_{\Omega}H(u_0)\phi_1\,dx}
{\int_{\Omega}f\phi_1\,dx},
$$
(3.11) follows.
\enddemo

\enddocument